\documentclass[letterpaper,12pt]{article}
\usepackage{amssymb}
\usepackage{amsmath}
\usepackage{makeidx}
\usepackage{amscd}

\def\vfl#1#2#3{\llap{$\scriptstyle #1$}
\left\downarrow\vbox to#3{}\right.\rlap{$\scriptstyle #2$}}

\newtheorem{theo}{Theorem}[section]
\newtheorem{prop}[theo]{Proposition}
\newtheorem{lem}[theo]{Lemma}
\newtheorem{cor}[theo]{Corollary}
\newtheorem{defi}[theo]{Definition}

\def \la {{\lambda}}

\def \si {{\sigma}}
\def \Ga {{\Gamma}}

\def \Pic {{\rm {Pic\,}}} 
\def \Div {{\rm {Div\,}}}

\def \A{{\mathbb A}}
\def \P{{\mathbb P}} 

\def \dim {{\rm{dim\,}}}

\def \Hom {{\rm {Hom}}}

\def \Pic {{\rm {Pic}}}
\def \GL {{\rm {GL}}}
\def \SL {{\rm {SL}}}

\def\ov{\overline}

\def \AA {{\rm A}}

\def \DD {{\rm D}}
\def \EE {{\rm E}}
\def \RR {{\rm R}}
\def \PP {{\rm P}}

\def \U {{\cal U}}

\def\G{{\bf G}}

\def\T{{\cal T}}
\def\sZ{{\cal Z}} 
\def\sY{{\cal Y}}

\def\lra{\longrightarrow}

\def\H{{\rm H}}

\def\O{{\cal O}}

\def\L{{\cal L }}

\def\O{{\cal O}}
\def\L{{\cal L}}

\def\sS{{\cal S}}

\def\si{\sigma}

\def\Ga{\Gamma}

\def\exp{{\rm exp}}

\newcommand{\bthe}{\begin{theo}}
\newcommand{\ble}{\begin{lem}}
\newcommand{\bpr}{\begin{prop}}
\newcommand{\bco}{\begin{cor}}
\newcommand{\bde}{\begin{defi}}
\newcommand{\ethe}{\end{theo}}
\newcommand{\ele}{\end{lem}}
\newcommand{\epr}{\end{prop}}
\newcommand{\eco}{\end{cor}}
\newcommand{\ede}{\end{defi}}

\title{On the equations for universal torsors over del~Pezzo surfaces}
\author{Vera Serganova and Alexei Skorobogatov} 
\date{\today} 
\begin{document}
\baselineskip=15pt
\maketitle
\smallskip

\centerline{\it \`a Jean-Louis Colliot-Th\'el\`ene}

\medskip

\section*{Introduction}

Universal torsors were invented by Jean-Louis Colliot-Th\'el\`ene and 
Jean-Jacques Sansuc; for smooth projective varieties $X$
with $\H^1(X,\O)=0$ they play the role similar to that of $n$-coverings
of elliptic curves.
The foundations of the theory of descent on torsors
were laid in a series of notes in
Comptes Rendus de l'Acad\'emie des Sciences de Paris in the second half 
of the 1970's, and a detailed account
was published in \cite{CS}. The theory has strong number theoretic
applications if the torsors can be described by explicit equations, and 
if the resulting system of equations can be treated using some other methods,
whether algebraic or analytic. Such is the case for surfaces fibred into
conics over $\P^1$: the universal torsors over these surfaces
are closely related to complete intersections of quadrics 
of a rather special kind. 
To describe them we use the following terminology. If $Z\subset\A_k^m$
is a closed subset of an affine space with a coordinate system
over a field $k$, then
the variety obtained from $Z$ by multiplying coordinates by non-zero
numbers will be called a dilatation of $Z$.
If exactly $n$ geometric fibres of the conic bundle $X\to\P^1$ are singular,
then there is a non-degenerate quadric $Q\subset\A_k^{2n}$ 
such that the 
universal torsors over $X$ are stably birationally equivalent 
to the product of a complete intersection of 
$n-2$ dilatations of $Q$, and a Severi--Brauer variety
(see \cite{CS}, Thm. 2.6.1).
This description was key for a plethora of applications of
the descent theory to the Hasse
principle, weak approximation, zero-cycles, R-equivalence
and rationality problems (the Zariski conjecture), see, e.g.,
\cite{CSS} and \cite{BCCS}. A stumbling block for a similar
treatment of cubic and more general smooth
del Pezzo surfaces without a pencil of rational curves is, possibly,
the absence of a satisfactory presentation of their universal torsors.
Known descriptions of universal torsors over
diagonal cubic surfaces (\cite{CS}, 2.5, \cite{CSK}, 10) 
lack the simplicity and the symmetry manifest in the conic bundle case.

One way to look at a non-degenerate quadric in $\A^{2n}$ is to think of 
it as a homogeneous space of the semisimple Lie group $G$ associated with
the root system $\DD_n$ which naturally appears in connection
with conic bundles with $n$ singular fibres, see e.g. \cite{KT}.
Indeed, over an algebraically closed field we can identify $Q$
with the orbit of the highest
weight vector of the fundamental $2n$-dimensional representation $V$ of $G$.
Then an `essential part' of the torsor is the intersection of $n-2$ dilatations
of this homogeneous space by the elements of a maximal torus in $\GL(V)$.

The aim of this paper is to generalize this description from the case
of conic bundles to that of del Pezzo surfaces. (Recall that these
two families exhaust all minimal smooth projective rational surfaces,
according to the classification of Enriques--Manin--Iskovskih.)
We build on the results of our previous paper \cite{I}, where
we studied split del Pezzo surfaces, i.e. the case when 
the Galois action on the set of exceptional curves is trivial.
The main result of \cite{I} is a construction of 
an embedding of a universal torsor over a split del Pezzo
surface $X$ of degree $5$, $4$, $3$ or $2$
into the orbit of the highest weight vector of a fundamental
representation of the semisimple simply connected Lie group 
$G$ which has the same root 
system as $X$, i.e. $\AA_4$, $\DD_5$, $\EE_6$ or $\EE_7$. 
This orbit is the punctured affine cone over $G/P$,
where $P\subset G$ is a maximal parabolic subgroup.
The embedding is equivariant with respect to the action of
the N\'eron--Severi torus $T$ of $X$, identified with a split maximal torus
of $G$ extended by $\G_m$.
In Theorem \ref{1.4} we describe universal torsors over 
split del Pezzo surfaces of degree $d$ as intersections of $6-d$
dilatations of the affine cone over $G/P$ by $k$-points of
a maximal torus of $\GL(V)$ which is
the centralizer of $T\subset\GL(V)$.
This gives a more conceptual approach to the equations
appeared previously in the work of Popov \cite{P} and Derenthal \cite{D}.
This approach can be called a {\it global description} of torsors
compared to their {\it local description} obtained by
Colliot-Th\'el\`ene and Sansuc  in \cite{CS}, 2.3.

For a general del Pezzo surface $X$
of degree $4$, $3$ or $2$ with a {\it rational point} we construct
an embedding of a universal torsor over $X$ into the same
homogeneous space as in the split case, but this time equivariantly
with respect to the action of a (possibly, non-split) maximal torus of $G$,
see Theorem \ref{nonsplit}. 
The case of del Pezzo surfaces of degree 5,
where a rational point comes for free by a theorem of Enriques and
Swinnerton-Dyer,
was already known (\cite{Sk}, Thm. 3.1.4).
The proof of Theorem \ref{nonsplit}
uses a recent result of Philippe Gille \cite{Gille} and 
M.S. Raghunathan \cite{R} which describes possible Galois actions 
on the character group of a maximal torus in a quasi-split algebraic group.
This result implies that the N\'eron--Severi torus $T$ of $X$
embeds into the same split group $G$ extended by $\G_m$, exactly as in the
case of a split del Pezzo surface.

 The condition on the
existence of a rational point on $X$ is not a restriction 
in the case of degree 5, but is clearly a restriction for smaller degrees,
limiting the scope of possible applications.
However, if $X$ is a del Pezzo surface of degree $4$, 
this condition is necessary as well as sufficient for our construction: 
if $X$ can be realized inside
a twisted form of the quotient of $G/P$ by a maximal torus,
then $X$ has a rational point, see Corollary \ref{c2} (i). 
Finally, in Corollary \ref{c2} (ii) we show that any del Pezzo surface
of degree $4$ with a $k$-point has a universal torsor which is
a dense open subset of the intersection of 
the affine cone over $G/P$ with its dilatation by a $k$-point of
the centralizer of $T$ in $\GL(V)$.

We recall the construction of \cite{I} in Section 1 alongside with
all necessary notation. In Section 2 we describe torsors over split
del Pezzo surfaces as intersections of dilatations of
the affine cone over $G/P$. In 
Section 3 we prove a uniqueness property used
in the proof of the main results in the non-split case in Section 4.

The ideas developed in this paper originate in the second author's
discussions with Victor Batyrev, to whom we are deeply grateful.

\section{Review of the split case}

\noindent{\it Preliminary remarks}\, Let $k$ be a field
of characteristic $0$ with an algebraic closure $\ov k$.

Let $V$ be a vector space over $k$, and let $T\subset \GL(V)$ be a split
torus, i.e. $T\simeq \G_m^n$ for some $n$.
Let $\Lambda\subset\hat T$ be the set of weights of $T$ in $V$, and
let $V_\la\subset V$ be the subspace of weight $\la$.
We have $V=\oplus_{\la\in\Lambda} V_\la$. Let $S$
be the centralizer of $T$ in $\GL(V)$, i.e.
$$S=\prod_{\la\in\Lambda} \GL(V_\la)\subset \GL(V).$$
In what follows we always assume that 
$\dim V_\la=1$ for all $\lambda\in\Lambda$; 
then $S$ is a maximal torus in $\GL(V)$.
Let $\pi_\la:V\to V_\la$ be the natural projection.
For $A\subset V$ we write $A^\times$ for the set of points
of $A$ outside $\cup\, \pi_\la^{-1}(0)$.

Let $r=4,\,5,\,6$ or $7$.
A split del Pezzo surface $X$ of degree $d=9-r$ is the blowing-up
of $\P^2$ in $r$ $k$-points in general position (i.e., no three
points are on a line and no six are on a conic). 
The Picard group $\Pic X$ is a free abelian group of rank $r+1$,
generated by the classes of exceptional curves on $X$.
Let $T=\G_m^{r+1}$. Once
an isomorphism $\hat T\tilde\lra\Pic X$ is fixed, $T$ is called
the N\'eron--Severi torus of $X$. 
A universal torsor $f:\T\to X$ is an $X$-torsor
with structure group $T$, whose type is the isomorphism
$\hat T\simeq\Pic X$ (see \cite{Sk}, p. 25).
We call a divisor in $\T$ an {\it exceptional divisor} if it is
the inverse image of an exceptional curve in $X$.

Now suppose that $\dim V$ equals the number of exceptional curves
on $X$. We can make an obvious but useful observation.

\ble \label{1.1}
Let $\T\to X$ be a universal torsor over a split del Pezzo surface $X$.
Let $\phi$ and $\psi$ be $T$-equivariant embeddings $\T\to V$
such that for each weight $\la\in\Lambda$
the divisors of functions $\pi_\la\phi$ and $\pi_\la\psi$ are equal
to the same exceptional divisor with multiplicity $1$.
Then $\psi=s\circ\phi$ for some $s\in S(k)$.
\ele
{\em Proof} 
Since $\T$ is a universal torsor we have $k[\T]^*=k^*$, hence
two regular functions with equal divisors differ by a non-zero
multiplicative constant. QED

\medskip

\noindent{\it Construction in the split case}\,
Let the pair consisting of a root system $\RR$ of rank $r$ 
and a simple root 
$\alpha$ be one of the pairs in the list
\begin{equation}
(\AA_4,\alpha_3), \quad (\DD_5,\alpha_5), \quad
(\EE_6,\alpha_6), \quad (\EE_7,\alpha_7).\label{list}
\end{equation}
Here and elsewhere in this paper we enumerate roots as in \cite{B}.
Let $G$ be the split simply connected simple group with split
maximal torus $H$ and root system $\RR$. 
Let $\omega$ be the fundamental weight dual to $\alpha$, and
let $V=V(\omega)$ be the irreducible $G$-module with the
highest weight $\omega$. It is known that $V$ is faithful and minuscule,
see \cite{B}. Let $P\subset G$ be the maximal parabolic 
subgroup such that $G/P\subset\P(V)$ is the orbit of the 
highest weight vector. The affine cone over $G/P$ is denoted by
$(G/P)_a$. 

It is easy to check that the $G$-module $S^2(V)$ is the direct sum
of two irreducible submodules $V(\omega_1)\oplus V(2\omega)$. 
For $r\leq 6$, $V(\omega_1)$ is a non-trivial irreducible $G$-module
of least dimension; it is a minuscule representation of $G$.
If $r=7$, then $V(\omega_1)$ is the adjoint representation;
it is quasi-minuscule, that is, all the non-zero weights 
have multiplicity 1 and form one orbit of $W$.
If $pr$ is the natural projection $V\to V(\omega_1)$, and
${\rm Ver}:V\to S^2(V)$ is the Veronese map $x\mapsto x^2$, then 
it is well known that
$(G/P)_a=(pr\circ {\rm Ver})^{-1}(0)$ (see \cite{BP}, Prop. 4.2 and references
there).

\medskip

Since the eigenspaces of $H$ in $V$ are 1-dimensional, 
$V$ has a natural coordinate
system with respect to which $S$ is the `diagonal' 
torus. 
Let the torus $T\subset S$ be the extension of $H$
by the scalar matrices $\G_m\subset\GL(V)$. Note that an eigenspace
of $H$ in $V$ is also an eigenspace of $T$, so that there is a natural
bijection between the corresponding sets of characters.

Let $V^{\rm sf}$ be the dense open subset of $V$ consisting 
of the points 
whose $H$-orbits are closed and whose stabilizers in $T$
are trivial. Let $(G/P)_a^{\rm sf}=(G/P)_a\cap V^{\rm sf}$.
In \cite{I} we constructed a $T$-equivariant closed embedding of 
$\T$ into $(G/P)_a^{\rm sf}$ 
such that each weight hyperplane section $\T\cap \pi_\la^{-1}(0)$
is an exceptional divisor with multiplicity 1. 
Then $X^\times=f(\T^\times)$ is the complement
to the union of exceptional curves on $X$.

We need to recall the details of this construction.
It starts with the case $(\RR,\alpha)=(\AA_4,\alpha_3)$ 
where the torsor $\T$ is the set of
stable points of $(G/P)_a$ which is the affine cone over
the Grassmannian ${\rm Gr}(2,5)$, see \cite{Sk}, 3.1. 
Thus $\T$ is open and dense in $(G/P)_a$ in this case.
As in \cite{I} we use dashes to denote the previous
pair in (\ref{list});
the previous pair of $(\AA_4,\alpha_3)$ is 
$(\AA_1\times\AA_2,\alpha^{(1)}_1+\alpha^{(2)}_2)$, 
though it will not be used.
For $r\geq 5$ we assume that a torsor $\T'\subset (G'/P')_a^{\rm sf}$
over a split del Pezzo surface of degree $10-r$
is already constructed, and proceed to construct $\T$
as follows.

Let $\Lambda_n\subset\Lambda$ be the set of weights $\la$ such that
$n$ is the coefficient of $\alpha$ in the decomposition
of $\omega-\la$ into a linear combination of simple roots.
Let $V_n=\oplus_{\la\in \Lambda_n}V_\la$, then
\begin{equation}
V=\bigoplus_{n\geq 0}V_n. \label{sum}
\end{equation}
The subspaces $V_n$ are $G'$-invariant. In fact, $V_n=0$ for $n>3$
so that $$V=V_0\oplus V_1\oplus V_2\oplus V_3,$$ 
and $V_3=0$ unless $r=7$. The degree 0 component
$V_0\simeq k$ is the highest weight subspace, and the degree 1 
component $V_1$ is isomorphic to $V'$ as a $G'$-module. 
The $G'$-module $V_2$ is irreducible with highest weight $\omega_1$.
For $r=7$ we have $V_3\simeq k$.

Recall from \cite{I} that $g_t=(t,1,t^{-1},t^{-2})$ is an element
of $T$, for any $t\in \ov k^*$.
Let $U\subset (G/P)_a$
be the set of points of $(G/P)_a$ outside
$(V_0\oplus V_1)\cup(V_2\oplus V_3)$. The natural
projection $\pi:V\to V_1$
defines a morphism $U\to V_1\setminus\{0\}$ which is the composition
of a torsor under $\G_m=\{g_t|t\in \ov k^*\}$ 
and the morphism inverse to the 
blowing-up of $(G'/P')_a\setminus\{0\}$ in 
$V_1\setminus\{0\}$ (\cite{I}, Cor. 4.2).
There is a $G'$-equivariant affine morphism
$\exp:V_1\to (G/P)_a$ such that $\pi\circ\exp={\rm id}$, and 
the affine cone over $\exp(V_1)$ is dense in $(G/P)_a$.
As in \cite{I} we write $\exp(x)=(1,x,p(x),q(x))$.
The map $p$ can be identified, up to a non-zero constant,
with the composition of the second Veronese map and the 
natural projection $S^2(V_1)\to V_2$, so that $(G'/P')_a=p^{-1}(0)$.
Since $V_2$ is the direct sum
of 1-dimensional weight spaces, it has a natural coordinate
system. The weight coordinates of $p(x)$ will be written as
$p_\mu(x)$, where $\mu\in W\omega_1$.

The choice of a point in $V^\times$ defines an isomorphism
$V^\times\simeq S$ compatible with the action of $S$. Using
this isomorphism we define a multiplication on $V^\times$,
and then extend it to $V$. However, none of our formulae will
depend of this isomorphism.
 
Suppose that $\T'\subset (G'/P')_a\subset V'=V_1$ is such that 
$f':\T'\to X'=\T'/T'$ is 
a universal torsor over a del Pezzo
surface $X'$, moreover, the $T'$-invariant hyperplane sections of $\T'$
are the exceptional divisors.
In \cite{I} we proved that for any $\ov k$-point $x_0$ in $\T'^\times$
there exists a non-empty open subset $\Omega(x_0)\subset (G'/P')_a^\times$,
whose definition is recalled in the beginning of the next section,
such that for any $y_0$ in $\Omega(x_0)$ the orbit $T'y_0$ is the scheme-theoretic
intersection $x_0^{-1}y_0\T'\cap (G'/P')_a$ (see Cor. 6.5 of \cite{I}).
Therefore, if $\T$ is the proper transform of $x_0^{-1}y_0\T'$ in $U$, 
then $X=\T/T$ is the blowing-up of $X'$ at the image of $x_0$. 
Consequently one proves that $\T\subset (G/P)_a^{\rm sf}$.
Equivalently, $\T$ can be defined as the affine cone (without zero)
over the Zariski closure of $\exp(x_0^{-1}y_0\T'\setminus T'y_0)$ 
in $(G/P)_a^{\rm sf}$.

The construction of an embedding of 
a universal torsor over $X$ into $(G/P)_a^{\rm sf}$
is a main result of \cite{I} (Thm. 6.1). The following corollary
to this theorem complements it by showing that our embedding is in 
a sense unique.

\bco \label{c1.2}
Let $\T\subset V^{\rm sf}$ be a closed $T$-invariant subvariety
such that $\T/T$ is a split del Pezzo surface and the weight
hyperplane sections of $\T$ are exceptional divisors with multiplicity $1$.
Then for some $s\in S(k)$ the torsor $s\T$ is a subset of $(G/P)_a$ 
obtained by our construction
(for some choice of a basis of simple roots of our root system $\RR$).
\eco
{\em Proof} The construction of \cite{I} recalled above produces
a universal torsor $\tilde\T$ over the same split del Pezzo surface $X$
inside $(G/P)_a$, satisfying the condition that the weight
hyperplane sections are the exceptional divisors with multiplicity $1$. 
The identifications of the exceptional curves on $X$ with the weights
of $V$ coming from $\T$ and $\tilde\T$ may be different, however
the permutation that links them is
an automorphism of the incidence graph of the
exceptional curves on $X$. It is well known (see \cite{M}) that the
group of automorphisms of this graph is the Weyl group ${W}$ of $\RR$.
Thus replacing $\tilde\T$ by its image under the action of
an appropriate element of ${W}$ (that is, a representative of
this element in the normalizer of $H$ in $G$), 
we ensure that the identification of the weights
with the exceptional curves is the same for both embeddings.
(The choice of this element in $W$ is equivalent to the choice of 
a basis of simple roots in our construction.)
The multiplicity 1 condition in the construction of \cite{I}
is easily checked by induction
from the case $r=4$ where we consider the Pl\"ucker coordinate
hyperplane sections of ${\rm Gr}(2,5)$.
It remains to apply Lemma \ref{1.1}. QED

\medskip

Let us recall some more notation.
For $\mu\in W\omega_1\subset \hat H$ we write $S^2_\mu(V)$
for the $H$-eigenspace of $S^2(V)$ of weight $\mu$, and $S^2_\mu(V)^*$
for the dual space. Let ${\rm Ver}_\mu$ 
be the Veronese map $V\to S^2(V)$
followed by the projection to $S^2_\mu(V)$.
For $r=6$ we write $S^3_0(V)$ for the zero weight
$H$-eigenspace in $S^3(V)$, and ${\rm Ver}_0:V\to S^3_0(V)$
for the corresponding natural map.

As in the previous corollary, we denote by $\T\subset V^{\rm sf}$ 
a closed $T$-invariant subvariety
such that $\T/T$ is a split del Pezzo surface and the weight
hyperplane sections of $\T$ are exceptional divisors with multiplicity $1$.
Let $I\subset k[V^*]$ be the ideal of $\T$,
$I_\mu=I\cap S^2_\mu(V)^*$, and, for $r=6$, let
$I_0=I\cap S^3_0(V)^*$. 

Let $\tilde \mu$ be the character by which $T$ acts on $S^2_\mu(V)$.
The $T$-invariant hypersurface in $\T$ cut by the 
zeros of a form from $S^2_\mu(V)^*\setminus I_\mu$ is mapped by 
$f:\T\to X$ to
a conic on $X$. The class of this conic in $\Pic X$, up to sign, 
is $\tilde\mu\in \hat T$ under the isomorphism $\hat T\simeq \Pic X$
given by the type of the torsor $f:\T\to X$ (see the comments
before Prop. 6.2 in \cite{I}).
All conics on $X$ in a given class are obtained in this way;
they form a 2-dimensional linear system,
hence the codimension of $I_\mu$ in $S^2_\mu(V)^*$ is 2.
Let $I_\mu^\perp\subset S^2_\mu(V)$ be the $2$-dimensional zero set of $I_\mu$.
The corresponding projective system defines a morphism
$f_\mu:X\to\P^1=\P(I_\mu^\perp)$ whose fibres are the conics 
of the class $\tilde\mu$.
The link between ${\rm Ver}_\mu$ and $f_\mu$ is described in the
following commutative diagram:
\begin{equation}
\begin{array}{ccccccc}V&\supset&(G/P)_a&\supset&\T&\to&X\\ 
\vfl{{\rm Ver}_\mu}{}{4mm}&&\vfl{}{}{4mm}&&\vfl{}{}{4mm}&&\vfl{f_\mu}{}{4mm}\\
S^2_\mu(V)&\supset&p_\mu^\perp&\supset&I_\mu^\perp\setminus \{0\}&\to&\P^1
\end{array} \label{d}
\end{equation}
Here $p_\mu^\perp$ is the zero set of $p_\mu\in S^2_\mu(V)^*$.

\ble \label{ld}
For $\mu\in W\omega_1$
the vertical maps in $(\ref{d})$ are surjective, and 
$\dim S^2_\mu(V)=r-1$.
\ele
{\em Proof} For the two right hand maps the statement is clear.
The map $V\to S^2_\mu(V)$ is surjective because all eigenspaces
of $T$ in $V$ are 1-dimensional. Since $\dim S^2_\mu(V)$
does not change if we replace $\mu$ by $w\mu$ for any $w\in W$,
to calculate $\dim S^2_\mu(V)$ we can assume that $\mu=\omega_1$.
But $\omega_1$ is a weight of $H'$ in $V_2$, so we have
$S^2_{\omega_1}(V)=S^2_{\omega_1}(V_1)\oplus (V_2)_{\omega_1}$, where
$\dim (V_{2})_{\omega_1}=1$. 
Starting with the case of the Pl\"ucker coordinates for $r=4$, one shows
by induction that $\dim S^2_{\omega_1}(V)=r-1$.
Hence $\dim S^2_\mu(V)=r-1$ and so $\dim p_\mu^\perp=r-2$ for any $\mu$.

To compute ${\rm Ver}_\mu((G/P)_a)$ we can continue to assume
that $\mu=\omega_1$.
If $x\in V_1$, then ${\rm Ver}_{\omega_1}$ sends
$\exp(x)\in(G/P)_a$ to
${\rm Ver}_{\omega_1}(x)+p_{\omega_1}(x)$. Thus
the projection of ${\rm Ver}_{\omega_1}((G/P)_a)$ to 
$S^2_{\omega_1}(V_1)={\rm Ver}_{\omega_1}(V_1)$, which is
a vector space of dimension $r-2$, is surjective. 
Hence ${\rm Ver}_{\omega_1}$ maps $(G/P)_a$ surjectively onto
$p_{\omega_1}^\perp$. QED

\medskip

Arguing by induction as in this proof,
it is easy to show starting with the case of 
the Pl\"ucker coordinates on ${\rm Gr}(2,5)$, that $p_\mu(x)$
is the sum of all the monomials of weight $\mu$ with non-zero
coefficients.

\section{Torsors over split del Pezzo surfaces}

Unless stated otherwise we assume that $r\geq 5$, 
so that $G'$ is of type $\AA_4$, $\DD_5$ or $\EE_6$. Recall that
we use dashes to denote objects related to the `previous'
root system.

Let $x_0$ be a $k$-point of $\T'^\times$. We define the dense open subset
$\Omega(x_0)\subset (G'/P')_a^\times$ as 
the set of $\ov k$-points $x$ such that $\exp(x_0^{-1}x\T'^\times)$
is not contained in $V\setminus V^\times$, that is, in 
the union of weight hyperplanes of $V$.
For $r=5$ or $6$ the set $\Omega(x_0)$
is the complement to the union of the closed subsets 
$$Z_\mu(x_0)=\{x\in (G'/P')_a\,|\,p_\mu(x_0^{-1}xu)\in I'_\mu\}$$
for all weights $\mu$ of $V_2$;
for $r=7$ one also removes the closed subset
$$Z_0(x_0)=\{x\in (G'/P')_a\,|\,q(x_0^{-1}xu)\in I'_0\}.$$
The condition $y_0\in \Omega(x_0)$ implies that for all $\mu$ the 
vectors ${\rm Ver}_\mu(x_0)$ and ${\rm Ver}_\mu(y_0)$ are not proportional.
Since $\dim S^2_\mu(V_1)^* =2+\dim I'_\mu$, we see that
for any $y_0\in\T'\cap \Omega(x_0)$ the subspace 
$I'_\mu\subset S^2_\mu(V_1)^*$ consists of the forms 
vanishing at $x_0$ and $y_0$.
The ideal $I'$ is generated by the $I'_\mu$, so
$\T'$ is uniquely determined by any two of its points
satisfying a certain open condition.

Recall that $\pi:V\to V_1$ is the natural projection, cf. (\ref{sum}).

\ble \label{1.0}
Let $x_0$ be a $k$-point of $\T'^\times$, 
let $y_0$ be a $k$-point of $\Omega(x_0)\cap\T'$, and let
$\T\subset (G/P)_a^{\rm sf}$ be the torsor 
defined by the triple $(\T',x_0, y_0)$ as described above.
Then we have the following statements.

{\rm (i)} The closed set $Z_\mu(x_0)\subset (G'/P')_a$ consists of
the points $x\in (G'/P')_a(\ov k)$ such that
$p_\mu(x_0^{-1}y_0x)=0$. For $r=7$ the closed set $Z_0(x_0)$ 
consists of the points $x\in (G'/P')_a(\ov k)$ such that
$q(x_0^{-1}y_0x)=0$.

{\rm (ii)} The open set $\Omega(x_0)\cap\T'$ is the inverse
image of the complement to all exceptional curves on $X'$
and to all conics on $X'$ passing through $f'(x_0)$. For $r=7$
one also removes from the cubic surface $X'\subset\P^3$
the nodal curve cut by the tangent plane to $X'$ at $f'(x_0)$.
We have $\T^\times=\pi^{-1}(\Omega(x_0)\cap\T')$.

{\rm (iii)} We have $t=\exp(x_0^{-1}y_0^2)\in \T^\times$.
\ele
{\em Proof} (i)
The inclusion of $Z_\mu(x_0)$ into the hypersurface given by
$p_\mu(x_0^{-1}y_0x)=0$ is clear: assigning the variable $u$
the value $y_0\in\T'$ we see that
$p_\mu(x_0^{-1}xu)\in I'_\mu$ implies that $p_\mu(x_0^{-1}x y_0)=0$.
Conversely, let us prove that every point $x$ of 
$(G'/P')_a$ satisfying the condition $p_\mu(x_0^{-1}y_0x)=0$,
is in $Z_\mu(x_0)$. Using Lemma \ref{ld} we see that the set of
quadratic forms $p_\mu(x_0^{-1}yu)$ on $V_1$
for a fixed $x_0$ and arbitrary $y\in(G'/P')_a$
is a vector subspace $L\subset S^2_\mu(V_1)^*$
of codimension 1, in fact this is the space of
forms vanishing at $x_0$. 
As was pointed out before the statement of the lemma,
$I'_\mu$ is the subspace of $L$ of codimension 1 consisting
of the forms vanishing at $y_0$. This proves the
desired inclusion.

Now let $r=7$. The inclusion of $Z_0(x_0)$ into the hypersurface 
$q(x_0^{-1}y_0x)=0$ is clear for the same reason as above.
Conversely, let $x\in (G'/P')_a(\ov k)$ be such that
$q(x_0^{-1}y_0x)=0$. We need to prove that $q(x_0^{-1}xu)$
vanishes for any $\ov k$-point $u$ of $\T'$.
In the end of the proof of Prop. 6.2 of \cite{I}
we showed that the dual space $\H^0(X',\O(-K_{X'}))^*$
is a 4-dimensional vector subspace of $S^3_0(V_1)$, so that
we have a commutative diagram similar to (\ref{d}):
$$\begin{array}{ccccc}
V_1&\supset&\T'&\to&X'\\
\vfl{{\rm Ver}_0}{}{4mm}&&\vfl{}{}{4mm}&&\vfl{\varphi}{}{4mm}\\
S^3_{0}(V_1)&\supset&\H^0(X',\O(-K_{X'}))^*\setminus\{0\}&\to& \P(\H^0(X',\O(-K_{X'}))^*)
\end{array}$$
where $\varphi$ is the anticanonical embedding $X'\hookrightarrow \P^3$.
In {\it loc. cit.} we also showed that for any $x\in(G'/P')_a(\ov k)$
the cubic form 
$q(x_0^{-1}xu)$, considered as a linear form on $S^3_0(V_1)$,
vanishes on the tangent space $T_{x_0}\simeq\P^2$ 
to $\varphi(X')\subset\P^3$ at $\varphi f'(x_0)$.
It is thus obvious that if $q(x_0^{-1}xu)$ vanishes at any point
of $\varphi(X')$ outside of $T_{x_0}$, then $q(x_0^{-1}xu)$
vanishes at any $\ov k$-point $u$ of $\T'$. But $\varphi f'(y_0)\notin T_{x_0}$,
otherwise $q(x_0^{-1}y_0z)=0$ for any $\ov k$-point $z$ of $(G'/P')_a$
contradicting the assumption that $y_0$ is in $\Omega(x_0)$.
Thus $q(x_0^{-1}xy_0)=0$ implies that $q(x_0^{-1}xu)\in I'_0$.

(ii) The geometric description of $\Omega(x_0)\cap\T'$ 
follows from \cite{I}, Cor. 6.3. Hence $\Omega(x_0)\cap\T'$
is obtained from $\T'$ by removing the images $\pi(E)$ of all
exceptional divisors $E\subset \T$, so that
$\pi(\T^\times)=\Omega(x_0)\cap\T'$.

(iii) Recall that $\exp(x)$ gives
a section of the natural morphism $\pi:\T\to x_0^{-1}y_0\T'$
over the complement to the fibre $T'y_0$. Thus $t\in \T$.
Since $y_0$ is in $\Omega(x_0)\cap\T'$ we see from (ii) that 
$t$ is in $\T^\times$. QED

\medskip

Let $\T\subset V^{\rm sf}$ be a closed $T$-invariant subvariety
such that $\T/T$ is a split del Pezzo surface and the weight
hyperplane sections of $\T$ are exceptional divisors with multiplicity $1$.
The torsor $\T$ defines an important subset of the torus $S$.
Namely, let $\sZ$ be the closed subset of $S$ consisting
of the points $s$ such that $s\T\subset (G/P)_a$. 
Equivalently, $\sZ=\bigcap_{x\in \T^\times(\ov k)}x^{-1}(G/P)_a^\times$.
The set $\sZ$ is $T$-invariant, since such are $(G/P)_a$ and $\T$.
In the case when $\T\subset (G/P)_a^{\rm sf}$, the variety $\sZ$
contains the identity element $1\in S(k)$. 

\ble \label{le}
Under the assumptions of Lemma $\ref{1.0}$ for
$r=4$ we have $\sZ=T$, and for $r\geq 5$ 
we have $\pi(\sZ)=y_0^{-1}\Omega(x_0)$ which is dense and open in 
$y_0^{-1}(G'/P')_a^\times$. 
The closed subvariety $\sZ\subset S$ is the affine cone (without zero) over
$t^{-1}\exp(x_0^{-1}y_0\Omega(x_0))$; in particular, $\sZ$
is geometrically integral, and $t^{-1}\T^\times\subset\sZ$.
For $r=5$ this inclusion is an equality.
\ele
{\em Proof} The statement in the case $r=4$ is clear since
$\T$ is dense in $(G/P)_a$, and the only elements of $S$ that leave
$G/P={\rm Gr}(2,5)$ invariant are the elements of $T$.
(Indeed, it is well known that the group of
relations among the classes of $10$ exceptional curves on $X$
is generated by the quadratic relations given by degenerate elements
of conic pencils on $X$. These quadratic relations are in a natural
bijection with the quadratic equations among the Pl\"ucker coordinates 
of ${\rm Gr}(2,5)$.)
Now assume that $r\geq 5$. For a fixed $x_0$, in order to construct an embedding
$\T\subset (G/P)_a$ we can choose any $y$ 
in the dense open subset $\Omega(x_0)\subset(G'/P')_a$. The embeddings 
defined by $(x_0,y_0)$ and $(x_0,y)$ satisfy the conditions of
Lemma \ref{1.1}. We obtain an
element $s\in \sZ$ such that $\pi(s)=y_0^{-1}y$. Thus $\pi(\sZ)$
contains $y_0^{-1}\Omega(x_0)$.

Let us prove that $\pi(\sZ)\subset y_0^{-1}(G'/P')_a^\times$.
Let $\pi_0:V\to V_0\simeq k$ be the natural projection.
Choose $y\in \T\subset (G/P)_a$ such that 
$\pi(y)=y_0\in \Omega(x_0)\subset (G'/P')_a^\times$.
By Lemma 4.1 of \cite{I} we have $\pi_0(y)=0$.
Thus $\pi_0(sy)=0$ for any $s\in \sZ$. But since $sy\in(G/P)_a$,
an inspection of cases in Lemma 4.1 of \cite{I} shows that
$\pi(sy)=\pi(s)y_0\in (G'/P')_a^\times$. 
Therefore, $\pi(\sZ)\subset y_0^{-1}(G'/P')_a^\times$.
Next, we note that $st\in (G/P)_a^\times$
(since $t\in\T^\times$ by Lemma \ref{1.0}). The coordinates
of the projection of $st$ to $V_2$ equal $p_\mu(\pi(s)x_0^{-1}y_0^2)$,
up to a non-zero constant, hence $p_\mu(\pi(s)x_0^{-1}y_0^2)\not=0$
for all $\mu$. But for $r\leq 6$ the open set $y_0^{-1}\Omega(x_0)
\subset y_0^{-1}(G'/P')_a^\times$ is given by 
$p_\mu(x_0^{-1}y_0^2u)\not=0$, by Lemma \ref{1.0} (i).
For $r=7$ a similar argument shows that $q(\pi(s)x_0^{-1}y_0^2)\not=0$.
Thus we obtain the equality $\pi(\sZ)=y_0^{-1}\Omega(x_0)$.

By Lemma \ref{1.0} (iii), $t=\exp(x_0^{-1}y_0^2)$ is in 
$\T^\times$ so we have $t\sZ\subset (G/P)_a^\times$.
Since $\sZ$ is invariant under the action of 
$\G_m=\{g_t|t\in \ov k^*\}$, we see from Lemma 4.1 of \cite{I} 
that $\sZ$ is a $\G_m$-torsor over $\pi(\sZ)=y_0^{-1}\Omega(x_0)$. 
Moreover, $t^{-1}\exp(x_0^{-1}y_0^2x)$ is a section of this torsor.
This proves that $\sZ$ is the affine cone over
$t^{-1}\exp(x_0^{-1}y_0\Omega(x_0))$.

If $r=5$, then $\Omega(x_0)$ is a dense open subset of $\T'$
as both sets are Zariski open in ${\rm Gr}(2,5)$. Thus the last
statement follows from Lemma \ref{1.0} (ii).
QED 

\medskip

This lemma implies that $\dim \sZ=2+\dim G'/P'$ which equals
$8,\,12,\,18$ for $r=5,\,6,\,7$, respectively.

\bde
$r-3$ points $z_0,\ldots,z_{r-4}$ in $\sZ(\ov k)$ are in general 
position if for any weight $\mu\in W\omega_1$ the vectors 
${\rm Ver}_\mu(z_i)$, $i=0,\ldots,r-4$, are linearly independent.
\ede

\ble \label{gp}
Let $\T\subset V^{\rm sf}$ be a closed $T$-invariant subvariety
such that $\T/T$ is a split del Pezzo surface and the weight
hyperplane sections of $\T$ are exceptional divisors with multiplicity $1$.
Then $\sZ$ contains $r-3$ $k$-points in general position.
More precisely, for any $k$-point $z_0$ of $\sZ$
the points $(z_1,\ldots,z_{r-4})\in \sZ(k)^{r-4}$
such that $z_0,z_1,\ldots,z_{r-4}$ are in general position,
form a dense open subset of $\sZ^{r-4}$.
\ele
{\em Proof} We first note that ${\rm Ver}_\mu(\sZ)$ is dense in a 
vector subspace of $S^2_\mu(V)$ of dimension $r-3$.
Indeed, assume without loss of generality that $\mu=\omega_1$.
Then, as in the proof of Lemma \ref{ld},
we have $S^2_{\omega_1}(V)=S^2_{\omega_1}(V_1)\oplus (V_2)_{\omega_1}$.
The image of $t\sZ$ consists of the points 
${\rm Ver}_{\omega_1}(x_0^{-1}y_0 u)+p_{\omega_1}(x_0^{-1}y_0 u)$, where
$u$ is in $\Omega(x_0)$, by Lemma \ref{le}. Since ${\rm Ver}_{\omega_1}$
sends $(G'/P')_a$ to a vector space of dimension $r-3$, by Lemma \ref{ld},
we see that ${\rm Ver}_{\omega_1}(t\sZ)$ is a dense subset of a vector
space of this dimension. Hence the same is true for $\sZ$.

We can choose the points
$z_1,\ldots,z_{r-4}$ in $\sZ(k)$ one by one, 
in such a way that $z_{n}$ is in the complement to
the union of the inverse images under ${\rm Ver}_\mu$
of the linear span of ${\rm Ver}_\mu(z_i)$, $i=0,\ldots,n-1$. 
This complement is non-empty since ${\rm Ver}_\mu(\sZ)$ 
is a Zariski dense subset of a vector space of dimension $r-3$. QED

\medskip

Equations for $\T$ have been given by Popov \cite{P} and Derenthal \cite{D}.
The following result gives a concise natural description of these 
equations, in terms of the well known equations of $(G/P)_a\subset V$.

\bthe \label{1.4}
Let $r=4,\,5,\,6$ or $7$. Every split del Pezzo surface $X$
of degree $9-r$ has a universal torsor $\T$ which is an
open subset of the intersection of $r-3$ dilatations of 
$(G/P)_a$ by $k$-points of the diagonal torus $S$.
In the above notation we have 
\begin{equation}
\T^\times=\bigcap_{z\in \sZ(\ov k)}z^{-1}(G/P)_a^\times
=\bigcap_{i=0}^{r-4}z_i^{-1}(G/P)_a^\times,
\label{**}
\end{equation}
where $z_0=1,\,z_1,\ldots,z_{r-4}$ are $k$-points of $\sZ$
in general position.
\ethe
{\em Proof} By Corollary \ref{c1.2} it is enough to prove the
theorem for $\T$ which satisfies the assumptions of Lemma \ref{1.0}.
The torsor $\T$ is clearly contained in 
the closed set $\sS=\cap_{s\in \sZ(\ov k)}s^{-1}(G/P)_a\subset V$.
Since $\T^\times$ is closed in $V^\times$, the density of $\T$ in $\sS$
implies $\T^\times=\sS^\times$. To prove this density
it is enough to show that $x_0^{-1}y_0\T'$ is dense in $\pi(\sS)$.

For $v\in V\otimes\ov k$ we write $v=(v_0,v_1,v_2,v_3)$, 
where $v_i\in V_i\otimes\ov k$.
Similarly, we write $s\in S(\ov k)$ as $(s_0,s_1,s_2,s_3)$, where 
$s_i\in \GL(V_i\otimes\ov k)$. In this notation the set
$$\bigcap_{s\in \sZ(\ov k)}\{(s_0^{-1}t,s_1^{-1}tx,s_2^{-1}tp(x),s_3^{-1}tq(x))|
x\in V_1\otimes\ov k,\ t\in \ov k^*\}$$
is dense in $\sS$. This set can also be written as
$$\bigcap_{s\in \sZ(\ov k)}\{(t,x,(ts_0)^{-1}s_2^{-1}p(s_1 x),
(ts_0)^{-2}s_3^{-1}q(s_1 x))|x\in V_1\otimes\ov k,\ t\in \ov k^*\}.$$
Since $(1,1,1,1)\in \sZ$, we see that
$\pi(\sS)$ is contained in the set of $x\in V_1\otimes\ov k$ 
such that for all $s\in \sZ(\ov k)$ we have 
$$s_0^{-1}s_2^{-1}p(s_1 x)=p(x).$$
Let $J\subset k[V_1^*]$ be the ideal
of $x_0^{-1}y_0\T'$, and $J_\mu=J\cap S^2_\mu(V_1)^*$. 
In the same way as $I'_\mu$, the ideal 
$J_\mu$ has codimension 2 in $S^2_\mu(V_1)^*$. Lemma \ref{ld}
implies that the linear span $L$
of the quadratic forms $p_\mu(yy_0^{-1}x)$ on $V_1$ for a fixed 
$y_0\in(G'/P')_a^\times$ and arbitrary $y\in(G'/P')_a^\times$
has codimension 1 in $S^2_\mu(V_1)^*$
(in fact, $L$ is the space of forms vanishing at $y_0$). 
Lemma \ref{le} implies that
$L$ coincides with the linear span of the quadratic forms 
$p_\mu(s_1 x)$, for all $s\in \sZ$. Hence the linear span of the forms 
$s_0^{-1}s_{2,\mu}^{-1}p_\mu(s_1 x)-p_\mu(x)$, for all $s\in \sZ$, has codimension
at most 2 in $S^2_\mu(V_1)^*$. However, the inclusion 
$x_0^{-1}y_0\T'\subset\pi(\sS)$ implies that this space is in $J_\mu$,
and thus coincides with $J_\mu$. This holds for every $\mu$, and the ideal
$J$ is generated by the $J_\mu$ (since the same is true for $I'$),
therefore $x_0^{-1}y_0\T'$ is dense in $\pi(\sS)$. 

Let us prove the second equality in (\ref{**}).
It is well known that the intersection of the ideal of 
$(G/P)_a$ with $S^2_\mu(V^*)$ is 1-dimensional;
let $P_\mu(u)$ be a non-zero element in this intersection.
Then $P_\mu(z_i u)$, where $z_0,\ldots,z_{r-4}$ are in general position,
span a vector space of dimension $r-3$ contained in the intersection of
the ideal of $\T$ with $S^2_\mu(V^*)$, which has the same dimension.
Thus $P_\mu(z_i u)$, $i=0,\ldots,r-4$, is a complete system of equations
of $\T$ of weight $\mu$. This completes the proof. QED

\medskip

\noindent{\it Remark} In the case $r=5$ the general position condition 
has a clear geometric interpretation. By the last claim of Lemma \ref{le}
we have $\T^\times=s\sZ$ for some $s\in S(k)$ well defined up to $T(k)$.
If $\T\subset (G/P)_a^{\rm sf}$, then $\sZ$ contains $1$, to that
$s$ is a $k$-point of $\T^\times$.
Then the previous theorem implies 
\begin{equation}
\T^\times=s\sZ=(G/P)_a^\times\cap r^{-1}s(G/P)_a^\times,
\label{*}
\end{equation}
where $r$ is a $k$-point in $\T^\times$ such that
$f(s)$ and $f(r)$ are points in $X^\times$ not contained in a conic
on $X$, cf. diagram (\ref{d}). Here $f(s)$ is uniquely
determined by $\T$, whereas $r$ can be any point
in the open subset of $X$ given by this condition.

\medskip

This remark can be seen as a particular case of the following
description of $\sZ$. For any $g$ and $h$ in $\T^\times(k)$ such that 
${\rm Ver}_\mu(h)$ and ${\rm Ver}_\mu(g)$ are not proportional
for any $\mu\in W\omega_1$,
we have $\sZ=g^{-1}(G/P)_a^\times\cap h^{-1}(G/P)_a^\times$.
The proof is similar to that of Theorem \ref{1.4}; we omit it here since
we shall not need this fact.

To construct $r-3$
points in $\sZ$ in general position is not hard, because the points of $\sZ$
are parameterized by polynomials. 
Indeed, decompose $V_1=V_{1,0}\oplus V_{1,1}\oplus V_{1,2}$
similarly to (\ref{sum}), and consider
the points $t^{-1}\exp(x_0^{-1}y_0\exp(v_i))$, where $v_1,\ldots,v_{r-3}$
in $V_{1,1}$ satisfy certain open conditions which are
easy to write down using Lemma \ref{le}.

\section{A uniqueness result}

The choice of $y_0$ plays the role of a `normalization'
for the embedding of a torsor into $(G/P)_a$. It is convenient
to choose these normalizations in a coherent way. Let $M_1,\ldots,M_r$
be $k$-points in general position in $\P^2$, and let $X_r$ be the blowing-up of $\P^2$ in $M_1,\ldots,M_r$.
The complement to the union of exceptional curves $X_r^\times\subset X_r$
can be identified with an open subset $\U\subset \P^2$.
Choose $u_0\in \U(k)$. At every step of our
inductive process we can choose
the points $y_0$ in the fibre of $\T'\to X'$ over $u_0$.
Thus we get a compatible
family of the $y_0$ (more precisely, of torus orbits)
that are mapped to each other by the surjective maps
$\T\to\T'$. In our previous notation, the point $t=\exp(x_0^{-1}y_0^2)$ 
must be taken for the point $y_0$ of the next step.

If $A$ is a subset of the torus $S$, then we denote by $\PP^n(A)\subset S$
the set of products of $n$ elements of $A$ in $S$.
We define $\PP^0(A)=T$.

\bpr \label{1.7}
Let $r$ and $n$ be integers satisfying $4\leq r\leq 7$, $0\leq n\leq r-4$.
Under the assumptions of Lemma $\ref{1.0}$, if at every step of our construction
we choose the points $y_0$ 
over a fixed point of $\U$, then we have the following statements:

{\rm (i)} $\PP^{n+1}(t^{-1}\T^\times)\subset t^{-1}(G/P)^\times_a$,

{\rm (ii)} $\PP^n(t^{-1}\T^\times)\subset \sZ$.
\epr
{\em Proof} (i) and (ii) are clearly equivalent. For $n=0$
the inclusion (i) is the main theorem of \cite{I}, 
and this also covers the case $r=4$.
Let $n\geq 1$. 
Recall that the projection $\pi$ maps $t^{-1}\T$ onto $y_0^{-1}\T'$.
Assume that we have the desired inclusions for $n-1$
and for both torsors $\T'$ and $\T$, namely
$$\PP^n(t^{-1}\T^\times)\subset t^{-1}(G/P)^\times_a, \quad
\PP^n(y_0^{-1}\T'^\times)\subset y_0^{-1}(G/P)^\times_a.$$
By Lemma 4.1 of \cite{I} every $\ov k$-point of $(G/P)^\times_a$
can be written as $g_x\cdot \exp(v)$, where $x\in \ov k^*$,
and $v\in V_1\otimes\ov k$. By the first inclusion in
induction assumption this is also true for
elements of $t\PP^n(t^{-1}\T^\times)$.
Since $\T$ is $g_x$-invariant, we have
$\exp(v)\in t\PP^n(t^{-1}\T^\times)$. 
On applying $\pi$ to both sides we deduce
$v\in x_0^{-1}y_0^2\PP^n(y_0^{-1}\T'^\times)$. 
Applying the second inclusion in induction
assumption we obtain $v\in x_0^{-1}y_0(G'/P')^\times_a$.
Therefore, 
$\PP^n(t^{-1}\T^\times)$ is contained in the affine cone over
$t^{-1}\exp(x_0^{-1}y_0(G'/P')^\times_a)$.
But this implies $\PP^n(t^{-1}\T^\times)\subset \sZ$, since 
$t\sZ$ is the intersection of the affine cone over
$\exp(x_0^{-1}y_0(G'/P')^\times_a)$ with $V^\times$, by 
the last statement of Lemma \ref{le}.
This proves (ii), and hence also (i).
QED

\bpr \label{1.8}
Let $r=4,\,5,\,6$ or $7$, and let $\T\subset V^{\rm sf}$ be a closed 
$T$-invariant subvariety
such that $X=\T/T$ is a split del Pezzo surface of degree $9-r$, 
and the weight
hyperplane sections of $\T$ are exceptional divisors with multiplicity $1$.
Let $\sZ\subset S$ be the closed subset
of points $z$ such that $z\T\subset (G/P)_a$.
Then there is a unique $s\in S(k)$ 
defined up to an element of $T(k)$,
such that $\PP^{r-4}(\T^\times)\subset s\sZ$.
\epr
{\em Proof} By Corollary \ref{c1.2}, 
up to translating $\T$ by an element of $S(k)$, we
can assume that $\T\subset(G/P)_a$ is obtained by our construction.
Thus the existence of $s$ follows from Proposition \ref{1.7}.
We prove the uniqueness by induction in $r$. For $r=4$ the
statement is clear, since the only elements of $S$ that leave
${\rm Gr}(2,5)$ invariant are the elements of $T$
(see the proof of Lemma \ref{le}).

Assume $r\geq 5$. 
By Lemma \ref{le}, $\PP^{r-4}(t^{-1}\T^\times)\subset s\sZ$ implies
$\PP^{r-4}(y_0^{-1}\T'^\times)\subset \pi(s)y_0^{-1}(G'/P')_a^\times$,
from which it follows that
$\PP^{r-5}(y_0^{-1}\T'^\times)\subset \pi(s)\sZ'$.
By induction assumption $\pi(s)$
is unique up to an element of $T'(k)$. Therefore, $s$
is unique up to an element of $T(k)$. QED

\medskip

\noindent{\it Remark} 
For $r=5$ the inclusion $\PP^{r-4}(\T^\times)\subset s\sZ$
is an equality by the last claim of Lemma \ref{le},
but this is no longer so for $r=6$ or $7$, for 
dimension reasons.

\section{Non-split del Pezzo surfaces}

Let $\Ga={\rm Gal}(\ov k/k)$.
Let $G$ be a split simply connected semisimple group over $k$
with a split maximal $k$-torus $H$ and the root system $\RR$.
There is a natural exact sequence of algebraic $k$-groups
\begin{equation}
1\to H\to N\to {W}\to 1, \label{N}
\end{equation}
where $N$ is the normalizer of $H$ in $G$, and 
${W}$ is the Weyl group of $\RR$. The action of $N$
by conjugation gives rise to an action of ${W}$ on the torus $H$.
Since $H$ is split, the Galois group
$\Ga$ acts trivially on ${W}$. Thus the continuous 1-cocycles
of $\Ga$ with values in $W$ are homomorphisms $\Ga\to W$, and
the elements of $\H^1(k,W)$ are homomorphisms $\Ga\to W$
considered up to conjugation in $W$.

\bthe[Gille--Raghunathan] \label{GR}
For any $\si\in\Hom(\Ga,{W})$ the
twisted torus $H_\si$ is isomorphic to a maximal torus of $G$.
\ethe
{\em Proof} See \cite{Gille}, Thm. 5.1 (b), or \cite{R}, Thm. 1.1. QED
\medskip

Recall from \cite{Serre}, I.5.4, that (\ref{N}) gives rise to 
the exact sequence of pointed sets
$$1\to N(k)\to G(k)\to (G/N)(k)
\buildrel{\varphi}\over{\hbox to 6mm{\rightarrowfill}} \H^1(k,N)\to \H^1(k,G).$$
(Note by the way that the last map here is surjective.)
The homogeneous space $G/N$ is the variety of maximal tori of $G$,
so that an equivalent form of the Gille--Raghunathan theorem 
is the surjectivity of the composite map 
$$(G/N)(k)\to \H^1(k,N)\to\H^1(k,{W})=\Hom(\Ga,{W})/{\rm Inn}({W}),$$
where ${\rm Inn}({W})$ is the group of inner automorphisms of ${W}$.
We fix an embedding of $H_\si$ as a maximal torus of $G$, this
produces a $k$-point $[H_\si]$ in $G/N$.
The choice of a $\ov k$-point $g_0$ in $G$
such that $g_0Hg_0^{-1}=H_\si$ defines a
1-cocycle $\rho:\Ga\to N(\ov k)$, $\rho(\gamma)=g_0^{-1}\cdot{{}^\gamma g_0}$,
which is a lifting of $\si\in Z^1(k,{W})=\Hom(\Ga,{W})$. 
We have $[\rho]=\varphi[H_\si]$ (\cite{Serre}, {\it ibidem}),
moreover, the image of $[\rho]$ in $\H^1(k,G)$ is trivial. 

\medskip

Let $G\to \GL(V)$ be an irreducible representation of $G$.
Define $T\subset \GL(V)$ as a torus
generated by $H$ and the scalar matrices $\G_m$.
The group $N$ acts by conjugation on $T$.
The twisted torus $T_\si$ is just the extension of $H_\si$ by scalar matrices.

Let $(G/P)_a\subset V$ be the orbit of the highest weight vector (with zero
added to it); $P\subset G$ is a parabolic subgroup, and $(G/P)_a$
is the affine cone over $G/P$.
The maximal torus $H_\si\subset G$ acts on $(G/P)_a$,
and so does $T_\si$. Define
$U_\si$ to be the dense open subset of $(G/P)_a$ consisting of the
points with closed $H_\si$-orbits and trivial stabilizers in $T_\si$.

The group $N\subset G$ acts on $V$ preserving $V^{\rm sf}$ 
and $V^\times$, 
thus giving rise to the action of ${W}$
on $V^{\rm sf}/T$ and on $V^\times/T$ by automorphisms 
of algebraic varieties
(not necessarily preserving some group structure on $V^\times/T$).
The action of $N$ preserves $(G/P)_a^{\rm sf}\subset V$, thus ${W}$
acts on $Y=(G/P)_a^{\rm sf}/T$. Hence we define the twisted forms
$(V^{\rm sf}/T)_\si$, $(V^\times/T)_\si$ and $Y_\si$.
The variety $(V^\times/T)_\si$
is an open subset of 
the quasi-projective toric variety $(V^{\rm sf}/T)_\si$, which contains
$Y_\si$ as a closed subset.

\ble \label{3.2}
The $k$-varieties $Y_\si$ and $U_\si/T_\si$ are isomorphic.
\ele
{\em Proof} Recall that $g_0\in G(\ov k)$ is a point such that
$\rho(\gamma)=g_0^{-1}\cdot{{}^\gamma g_0}\in Z^1(k,N)$
is a cocycle 
that lifts $\si\in Z^1(k,{W})=\Hom(\Ga,{W})$.
The group $N$ acts on the homogeneous space $(G/P)_a$ as a subgroup of $G$,
so we can define
$(G/P)_{a,\rho}$ as the twist of $(G/P)_{a}$ by $\rho$. 
It is immediate to check that
the map $x\mapsto g_0 x$ on $\ov k$-points of $(G/P)_{a}$
gives rise to an isomorphism of $k$-varieties 
$(G/P)_{a,\rho}\tilde\lra (G/P)_{a}$. If $G_\rho$
is the inner form of $G$ defined by $\rho$, then $G_\rho$ acts on
$(G/P)_{a,\rho}$ on the left. The embedding $H\hookrightarrow G$
gives rise to an embedding $H_\sigma\hookrightarrow G_\rho$,
so that $T_\sigma$ acts on $(G/P)_{a,\rho}$ on the left.
On the other hand, $T_\sigma$ also acts on $(G/P)_a$ on the left.
It is straightforward to check that
the isomorphism $(G/P)_{a,\rho}\tilde\lra (G/P)_{a}$ is
$T_\sigma$-equivariant.

Let $(G/P)_{a,\rho}^{\rm sf}$ be the subset of $(G/P)_{a,\rho}$
consisting of the
points with closed $H_\si$-orbits with trivial stabilizers in $T_\si$.
The closedness of orbits and the triviality
of stabilizers are conditions on $\ov k$-points, hence we obtain 
a $T_\si$-equivariant 
$k$-isomorphism $(G/P)_{a,\rho}^{\rm sf}\tilde\lra U_\si$. It
descends to an isomorphism $Y_\si\tilde\lra U_\si/T_\si$.
This proves the lemma. QED

\bco \label{c1}
For any homomorphism $\si:\Ga\to{W}$ the twisted
variety $Y_\si$ has a $k$-point, and so does $(V^\times/T)_\si$.
\eco
{\em Proof} Since $k$ is an infinite field, any dense open subset
of $(G/P)_a$ contains $k$-points.
Thus $Y_\si^\times(k)\not=\emptyset$, but this is a subset of 
$(V^\times/T)_\si$, so that this variety has a $k$-point. QED

\medskip

\noindent{\it Remark} This approach via the Gille--Raghunathan theorem
generalizes a key ingredient in
the second author's proof of the Enriques--Swinnerton-Dyer theorem
that every del Pezzo surface of degree 5 has a $k$-point, from 
quotients of Grassmannians by the action of a maximal torus
to quotients of arbitrary
homogeneous spaces of quasi-split semisimple groups. We plan to return
to this more general statement in another publication.

\bigskip

We now assume that $\RR$ is the root systems of rank $r$ in (\ref{list}),
and that the highest weight of the $G$-module $V$ is the fundamental
weight dual to the root indicated in (\ref{list}). Then $V$
is minuscule, so that the centralizer $S$ of $H$ in $\GL(V)$ is a torus.
Let $R=S/T$. We obtain an exact sequence of $k$-tori:
\begin{equation}
1\to T\to S\to R\to 1.
\end{equation}
The group $N$ acts by conjugation on $T$ and hence also on $S$ and $R$.
The connected component of 1 acts trivially, so we obtain
an action of ${W}$ on these tori (preserving the group structure). 
On twisting $T$, $S$ and $R$ by $\si$ we obtain
an exact sequence of $k$-tori:
\begin{equation}
0\to T_\si\to S_\si\to R_\si\to 0. \label{twist}
\end{equation}
Note in passing that the character group $\hat S$ has an obvious ${W}$-invariant basis, which
gives rise to a Galois invariant basis of $\hat S_\si$.
In other words, $S_\si$ is a quasi-trivial torus;
in particular, $\H^1(k,S_\si)=0$ as follows from Hilbert's theorem 90.
Note also that $V^\times/T$ is a torsor under $R$, so that $(V^\times/T)_\si$
is a torsor under $R_\si$. By Corollary \ref{c1} this torsor is
trivial, that is, there is a (non-canonical) isomorphism 
$(V^\times/T)_\si\simeq R_\si$.

\medskip

Let $X$ be a del Pezzo surface over $k$, not necessarily split,
of degree $9-r$, where $r$ is the rank of the root system $\RR$. 
Let $\ov X$ be the surface obtained from $X$ by extending 
the ground field from $k$ to $\ov k$. We write
$X^\times$ for the complement to the union of exceptional curves on $X$.
Our construction identifies $\hat S$ with the free abelian group
$\Div_{\ov X\setminus \ov {X^\times}}\ov X$
generated by the exceptional curves on $\ov X$, and
$\hat T$ with $\Pic\ov X$ 
(via the type of the universal torsor $\T\to X$).
The Galois group permutes the exceptional curves on $\ov X$, thus
defining a homomorphism 
$\si_X:\Ga\to{W}$, where ${W}$ is the Weyl group of $\RR$.
This homomorphism is well defined up to conjugation in ${W}$,
so we have a well defined class $[\si_X]\in\H^1(\Ga,{W})$,
where $\Ga$ acts trivially on ${W}$. 

We now assume $\si=\si_X$. Then we get isomorphisms of $\Ga$-modules
$$\hat S_\si=\Div_{\ov X\setminus \ov {X^\times}}\ov X, \quad \quad
\hat T_\si=\Pic\ov X,$$
thus $T_\si$ is the N\'eron--Severi torus of $X$.
The kernel of the obvious surjective map
$\Div_{\ov X\setminus \ov {X^\times}}\ov X\to \Pic\ov X$
is $\ov k[X^\times]^*/\ov k^*$, hence
the dual sequence of (\ref{twist}) coincides with the 
natural exact sequence of $\Ga$-modules
\begin{equation}
0\to \ov k[X^\times]^*/\ov k^*\to 
\Div_{\ov X\setminus \ov {X^\times}}\ov X\to \Pic\ov X\to 0. \label{pics}
\end{equation}
There is a natural bijection between the morphisms $X^\times\to R_\si$ 
and the homomorphisms of $\Ga$-modules $\hat R_\si\to \ov k[X^\times]^*$.
Universal torsors on $X$ exist if and
only if the exact sequence of $\Ga$-modules
\begin{equation}
1\to\ov k^*\to \ov k[X^\times]^*\to\ov k[X^\times]^*/\ov k^*\to 1 \label{ls}
\end{equation}
is split (\cite{Sk}, Cor. 2.3.10). Any splitting of this
sequence gives a map
$$\hat R_\si=\ov k[R_\si]^*/\ov k^*=\ov k[X^\times]^*/\ov k^*\to \ov k[X^\times]^*,$$
and hence defines a morphism $\phi:X^\times\to R_\si$. 
By the `local description of torsors' (see \cite{CS}, 2.3 
or \cite{Sk}, Thm. 4.3.1) the restriction
of a universal $X$-torsor to $X^\times$ is the pull-back
of the torsor $S_\si\to R_\si$ to $X^\times$ via $\phi$.
Moreover, this gives a bijection between the splittings of
(\ref{ls}) and the universal $X$-torsors.
In our case it is easy to see that $\phi$ is an embedding.
The isomorphism $\hat R_\si=\ov k[X^\times]^*/\ov k^*$
comes from our construction, thus after extending the ground field
to $\ov k$, the morphism $\phi$ coincides, up to translation
by a $\ov k$-point of $R$, with the embedding of
$\ov X^\times$ into $(V\otimes_k\ov k)^\times/ \ov T$ 
obtained from the embedding
$\ov\T^\times\subset\ov V^\times$.

\bthe \label{nonsplit}
Let $r=4,\,5,\,6$ or $7$.
Let $X$ be a del Pezzo surface of degree $9-r$
with a $k$-point, and let
$\si\in\H^1(\Ga,{W})$ be the class defined by the action of 
the Galois group on the exceptional curves of $X$. 
There exists an embedding $X\hookrightarrow Y_\si$
such that the divisors in $Y_\si\setminus Y_\si^\times$ cut the exceptional
curves on $X$ with multiplicity $1$.
The restriction of $U_\si\to Y_\si$ to $X\subset Y_\si$ is 
a universal $X$-torsor whose type is the isomorphism 
$\hat T_\si=\Pic\ov X$.
\ethe
{\em Proof} From Corollary \ref{c1} we get an embedding 
$Y_\si\hookrightarrow R_\si$, which becomes unique if we further
assume that a given $k$-point of $Y_\si$ goes to the identity element
of $R_\si$.

Since $X(k)\not=\emptyset$, there is a unique embedding
$\phi:X^\times\to R_\si$ such that the induced map
$\phi^*:\hat R_\si\to \ov k[X^\times]^*$ is a lifting of the
isomorphism $\hat R_\si=\ov k[R_\si]^*/\ov k=\ov k[X^\times]^*/\ov k^*$,
and $\phi$ sends a given $k$-point of $X^\times$ to $1$.

Let $\L$ be the $k$-subvariety of the torus $R_\si$ whose points
are $c\in R_\si(\ov k)$ such that $cX^\times\subset Y_\si^\times$, where
the multiplication is the group law of $R_\si$.
To prove the first statement we need to show that $\L(k)\not=\emptyset$.
Let $\PP^n(X^\times)$ be the $k$-subvariety of $R_\si$ whose
$\ov k$-points are products of $n$ elements
of $X^\times(\ov k)$ in $R_\si(\ov k)$. The surface
$\ov X$ is split,
hence it follows from Proposition \ref{1.8} that
there exists a unique $c_0\in R_\si(\ov k)$ such that
$\PP^{r-4}(X^\times)(\ov k)\subset c_0\L(\ov k)$.
But since $\PP^{r-4}(X^\times)$ and $\L$ are 
subvarieties of $R_\si$ defined over $k$
we conclude that $c_0$ is a $k$-point.
If $m$ is $k$-point of $X^\times$, 
then $c_0^{-1}m^{r-4}$ is a $k$-point of $\L$, as required. 

To check that the restriction of $U_\si\to Y_\si$ to $X\subset Y_\si$
is a universal torsor we can go over to $\ov k$ where it follows from
our main theorem in the split case. QED

\medskip

\noindent{\it Remark} Let $X$ be a del Pezzo surface with a $k$-point,
of degree 5, 4, 3 or 2. Although $(G/P)_a$ contains some universal $X$-torsor,
other universal $X$-torsors of the same type are naturally 
embedded into certain twists of $(G/P)_a$. Indeed,
all torsors of the same type are obtained from any of them by twisting
by the cocycles in $Z^1(k,T_\si)$. The natural map
$H_\si\to T_\si$ gives a surjection $\H^1(k,H_\si)\to\H^1(k,T_\si)$
since $\H^1(k,\G_m)=0$ by Hilbert's theorem 90. Therefore, if 
$\T\subset (G/P)_a$ and $\theta\in Z^1(k,H_\si)$, then 
$\T_\theta$ is contained in the twist of $(G/P)_a$ by the 1-cocycle 
in $Z^1(k,G)$ coming from $\theta\in Z^1(k,H_\si)$. By general theory
(\cite{Serre}, I.5)
this twist is a left homogeneous space of the inner
form of $G$ defined by $\theta$. (In the case of ${\rm Gr}(2,5)$
we only obtain twists that are isomorphic to ${\rm Gr}(2,5)$ because
$\H^1(k,\SL(n))=0$.) We note that by Steinberg's theorem
every class in $\H^1(k,G)$ comes from $\H^1(k,H_\si)$
for some maximal torus $H_\si\subset G$.

\medskip

\bco \label{c2}
Let $X$ be a del Pezzo surface of degree $4$
such that universal $X$-torsors exist.
Let $\si\in\H^1(\Ga,{W})$ be the class defined by the action of 
the Galois group on the exceptional curves of $X$. 

{\rm (i)}
\, $X^\times$ and $Y_\si^\times$ are $k$-subvarieties of $R_\si$.
Moreover, $Y_\si$ contains $cX$ for some $c\in R_\si(k)$
if and only if $X$ has a $k$-point.

{\rm (ii)} If $X\subset Y_\si$, then
$X=Y_\si\cap m Y_\si$ for some $m\in R_\si(k)$.
\eco
{\em Proof.} (i) In view of Theorem \ref{nonsplit} it remains to prove the
`only if' part. We note that $Y_\si$ embeds into $R_\si$ by 
Corollary \ref{c1}. The existence of universal $X$-torsors implies that
$X$ embeds into $R_\si$, as was discussed before Theorem \ref{nonsplit}.
For $r=5$ the inclusion $X^\times\subset c_0\L$ from the proof of Theorem
\ref{nonsplit} is an equality
by the remark in the end of Section 3. Hence 
if $\L$ has a $k$-point, then so does $X^\times$. 

(ii) A del Pezzo surface of degree $4$ with a $k$-point is known to 
be unirational (i.e., is dominated by a $k$-rational variety, see \cite{M}).
Thus $X$, and hence also $\L=c_0^{-1} X^\times$ 
contains a Zariski dense set of $k$-points. The variety $X$
is contained in $Y_\si\cap h^{-1}c_0Y_\si$ for any
$h\in X^\times(k)$, and this inclusion is an equality for any $h$
in a dense open subset of $X$, see the remark after the proof of Theorem
\ref{1.4}. QED

\medskip

\noindent{\it Remark} By the previous remark 
an arbitrary universal torsor over a del Pezzo surface $X$ 
of degree $4$
with a $k$-point embeds into a twisted form of $(G/P)_a$
by a cocycle coming from $\theta\in Z^1(k,H_\si)$. 
This twisted form $(G/P)_{a,\theta}$ is naturally a subset of a vector space
(non-canonically isomorphic to $V$) acted on by $S_\si$. Thus
any universal torsor over $X$ is an open subset of the intersection of 
two $k$-dilatations of $(G/P)_{a,\theta}$.

\noindent Department of Mathematics, University of California,
Berkeley, CA, 94720-3840 USA
\medskip

\noindent serganov@math.berkeley.edu

\bigskip

\noindent Department of Mathematics, South Kensington Campus, 
Imperial College London, SW7 2BZ England, U.K.

\smallskip

\noindent Institute for the Information Transmission Problems, 
Russian Academy of Sciences, 19 Bolshoi Karetnyi, 
Moscow, 127994 Russia
\medskip

\noindent a.skorobogatov@imperial.ac.uk

\end{document}